\documentclass[12pt]{amsart}

\usepackage{amssymb}
\usepackage{amsthm}
\usepackage{amsmath}
\usepackage{cite}
\usepackage{verbatim}

\newtheorem{theorem}{Theorem}

\newcommand{\maxprime}{83}
\newcommand{\maxcomp}{52}
\newcommand{\Z}{\mathbb{Z}}
\newcommand{\Q}{\mathbb{Q}}

\def\A{\mathbb A}

\def\R{\mathbb R}
\def\C{\mathbb C}
\def\Q{\mathbb Q}
\def\Z{\mathbb Z}

\def\L{\mathcal L}

\def\cpi{{}^c\pi}
\def\tpi{\widetilde\pi}
\def\Gal{\mbox{Gal}}
\def\As{\mbox{As}}
\def\cusp{\mbox{cusp}}

\def\diag{\mbox{diag}\,}

\def\det{\mbox{det}}

\def\BorelSerre{\text{BS}}
\def\Eis{\text{Eis}}
\def\cusp{\text{cusp}}
\def\Frob{\text{Frob}}
\def\Gal{\text{Gal}}
\def\IIa{\text{IIa}}
\def\IIb{\text{IIb}}
\def\IIIa{\text{IIIa}}
\def\IIIb{\text{IIIb}}
\def\IV{\text{IV}}

\theoremstyle{plain}
\newtheorem{conjecture}{Conjecture}

\theoremstyle{definition}
\newtheorem{example}{Example}

\begin{document}

\title{Cohomology of Congruence Subgroups of $SL(4,\Z)$ II}

\author{Avner Ash} \address{Boston College\\ Chestnut Hill, MA 02445}
\email{Avner.Ash@bc.edu} \author{ Paul E. Gunnells}
\address{University of Massachusetts\\ Amherst, MA 01003}
\email{gunnells@math.umass.edu} \author{Mark McConnell}
\address{WANDL, Inc.\\ 27 Wolf Hill Drive, Warren, New Jersey 07059}
\email{mmcconnell@wandl.com} \thanks{The first author wishes to thank
the National Science Foundation for support of this research through
NSF grants numbers DMS-0139287 and DMS-0455240.  The second author
wishes to thank the National Science Foundation for support of this
research through NSF grants numbers DMS-0245580 and DMS-0401525.  We
thank the mathematics departments of Columbia University and the
University of Massachusetts Amherst for computers used to perform
these computations.  We also thank Armand Brumer, Dinakar
Ramakrishnan, and Uwe Weselmann
for helpful conversations.
groups, Hecke operators, Eisenstein cohomology, Siegel modular forms.}

\date{31 May 2007}

\begin{abstract}
In a previous paper \cite{AGM} we computed cohomology groups $H^{5}
(\Gamma_{0} (N), \C)$, where $\Gamma_{0} (N)$ is a certain congruence
subgroup of $SL (4, \Z)$, for a range of levels $N$.  In this note we
update this earlier work by extending the range of levels and describe
cuspidal cohomology classes and additional boundary phenomena found
since the publication of \cite{AGM}.  The cuspidal cohomology classes
in this paper are the first cuspforms for $GL(4)$ concretely
constructed in terms of Betti cohomology.
\end{abstract}

\maketitle

\section{Introduction}\label{intro}

In this paper we extend the computations in \cite{AGM} of the
cohomology in degree 5 of congruence subgroups $\Gamma_0(N) \subset
SL(4,\Z)$ with trivial $\C$ coefficients to more levels, up to level
\maxprime.  We also compute Hecke operators on these cohomology groups
and use the Hecke eigenvalues to identify the cohomology eigenclasses
as either Eisenstein or cuspidal.  We remind the reader that
$\Gamma_0(N)$ is the subgroup of $SL(4,\Z)$ consisting of matrices
whose last row is congruent modulo $N$ to $(0,0,0,*)$.  We say that
$\Gamma_0(N)$ is ``modeled" on the $(3,1)$ parabolic subgroup of
$SL(4)/\Q$.  Also recall that a cohomological cuspidal automorphic
representation contributes to the cohomology of $\Gamma_0(N)$ exactly
in degrees 4 and 5.

The size of the matrices and the complexity of computing the Hecke
operators are greater the larger $N$ or the more composite $N$
is. Similarly the size of the computation of the Hecke operators at a
prime $\ell $ increases dramatically as a function of both $\ell$ and
$N$.  Therefore after $N=\maxcomp$ we stopped computing for composite
$N$ but were able to continue for prime $N$ up to level \maxprime.
Similarly the size of the computation of the Hecke operators at a
prime $\ell$ increase dramatically with $\ell$, so that in fact for
the new levels in this paper, we computed the Hecke operators only for
$\ell = 2$ and in a few cases for $\ell = 3,5$.

For levels $N = 2$ through $31$, we have checked our results by redoing
the computations but with coefficients in $\Z$.  In this way, we have also
identified non-trivial torsion classes in $H^5(\Gamma_0(N), \Z)$ for some
levels $N$.  These torsion classes and their relationship to Galois
representations will be studied in a future paper.

Working with $\Z$ coefficients is more difficult than with coefficients in a
finite field, because the size of the integers in the intermediate steps
of the calculations tends to grow exponentially. This is why we stopped
using $\Z$ at $N = 31$.  For higher levels, we worked over the finite
fields $\Z/31991\Z$ or $\Z/12379\Z$.

Unlike the earlier paper, our new results include cuspforms.  They
also confirm the observed patterns of Eisenstein liftings from the
cohomology of the Borel--Serre boundary of the locally symmetric
space for $\Gamma_0(N)$ as explained in \cite{AGM}.
We refer the reader to \cite{AGM} for a detailed explanation of why we
look in degree 5, the interaction with the cohomology of the
Borel--Serre boundary, and how we perform the computations.

\medskip \noindent \emph{Cuspforms.}  We discovered cuspforms at
levels 61, 73, and 79.  All these levels are prime.  To the best of
our knowledge, these are the first concretely constructed cuspforms
for $GL(4)$ in the sense of Betti cohomology.  Each cuspform appears
with multiplicity two in the cohomology, viewed as a module for the
Hecke operators, the eigenvalues being rational integers.  In Section
\ref{cusp} we explain this.  Theorem \ref{Dinakar} 
asserts that these cuspforms must be functorial
liftings from holomorphic Siegel modular forms of weight 3 on
$GSp(4)/\Q$.

It would be interesting to see a construction of these Siegel modular
forms.  There are several potential approaches.  One is to construct
them using theta series.  Such theta series would be on the congruence
subgroups modeled on the Klingen parabolic subgroup in $GSp(4)/\Q$.
Unfortunately all the work we know of that might have been relevant,
for example \cite{PY}, concerns congruence subgroups modeled on the
Siegel parabolic subgroup.  These latter theta series lead to Siegel
modular forms that can appear in the cohomology of the congruence
subgroups of $SL(4,\Z)$ modeled on the $(2,2)$ parabolic subgroup of
$SL(4)/\Q $.  We plan to investigate the cohomology of such subgroups
in future work.  Another possible construction of the desired Siegel
modular forms is by computing the cohomology of congruence subgroups
of $Sp (4, \Z)$, since holomorphic Siegel modular forms of weight $3$
will contribute to the cohomology of these groups.  Finally, one could
also try to isolate the motives corresponding to the cuspidal
automorphic forms we found.  These might be found either as factors in the
cohomology of the appropriate Siegel modular variety, or by other
means, as in the work of van Geemen and Top \cite{vgt}.

We remark that Ibukiyama \cite{ibuk} conjectured (and has recently
announced a proof of) a formula that describes the dimensions of
weight three cuspidal Siegel modular forms on the paramodular groups
of prime level.  The paramodular group of level $N$ is a congruence
subgroup of $Sp (4, \Q)$ that contains the congruence subgroup of
level $N$ based on the Klingen parabolic subgroup.  Gritsenko
\cite{grit} has constructed a lift from Jacobi forms to Siegel modular
forms on the paramodular group.  Brumer observed that the first levels
where the forms predicted by Ibukiyama are not accounted for by
Gritsenko's lifts are 61, 73, and 79, and in these cases the subspace
of lifts has codimension one.  Thus we expect that our classes will
prove to be concrete realizations of the lifts of these Siegel modular
forms.

It would be most interesting to discover cohomology classes in
$H^5(\Gamma_0(N),\C)$ corresponding to cuspforms that are not lifts
from any smaller group, but these have not shown up yet in our
computations.  Each such cuspidal Hecke eigenclass would give rise to
a 4-dimensional subspace of the cohomology over $\C$, namely two
subspaces each of multiplicity two, with Hecke eigenvalues the complex
conjugates of each other (see Section \ref{cusp}).

\medskip \noindent \emph{Eisenstein series.}
We continue to observe that all weight 2 newforms for $GL(2)/\Q$ of
level dividing $N$ lift as cohomology of the Borel--Serre boundary
into the cohomology of our $\Gamma_0(N)$.  Only certain weight 4
newforms $f$ were observed to lift.  We conjecture that such a form
lifts if and only if the sign in its functional equation is negative
(see Conjecture \ref{eisenstein.conjecture} below for more details).  The
connection between this sign and our observed lifting phenomenon was
pointed out to us by U.~Weselmann.

For prime levels $N$, we have also observed that cohomology classes
with trivial coefficients of level $N$ attached to cuspidal
automorphic representations of $GL (3)/\Q$ lift to $H^{5} (\Gamma_{0}
(N), \C )$, again via the cohomology of the Borel--Serre boundary.

When the level $N$ is a square, we observed that some cohomology for
minimal faces of the Borel--Serre boundary lifts to $H^{5}$.  The
details of this phenomenon are currently unclear.\label{minimalref}

Eisenstein cohomology, originally introduced by G.~Harder
(cf.~\cite{harder-icm}), has been investigated extensively by
J.~Franke, J.~Rohlfs, J.~Schwermer, B.~Speh, and others.  However, it
is a difficult open problem to compute in all detail the cohomology of
the Borel--Serre boundary for $GL(4)$ for general $N$.  Even if that
were done, current results in Eisenstein cohomology don't appear to be
fine enough even to check Conjecture \ref{eisenstein.conjecture},
which is only for prime level.  One problem is that 
trivial coefficients is much harder to handle than irreducible
coefficient modules with regular highest weights.

In the final section of the paper we provide some tables of results.
Table \ref{bettitab} shows the new Betti numbers we computed and
extends the data in \cite{AGM}.\footnote{Level 49 was incorrectly
reported in \cite{AGM}.  Levels 55, 67, 71 were conjectured in
\cite{AGM}.}  Table \ref{grtab} shows the Hecke polynomials for the
Eisenstein classes we found that are covered by Conjecture
\ref{eisenstein.conjecture}.  Finally Table \ref{cusptab} gives the
Hecke polynomials and eigenvalues for the cuspidal classes we found.

\section{Eisenstein classes}\label{eis}

Let $\xi \in H^{5} (\Gamma_{0} (N), \C)$ be a Hecke eigenclass.
Recall \cite[\S1.1]{AGM} that for us this means $\xi$ is an
eigenvector for certain operators 
\[
T (l,k) \colon H^{5} (\Gamma_{0}
(N), \C) \rightarrow H^{5} (\Gamma_{0} (N), \C),
\]
where $k=1,2,3$ and
$l$ is a prime not dividing $N$.  These operators correspond to the
double cosets $\Gamma_{0} (N) D (l,k) \Gamma_{0} (N)$, where $D (l,k)$
is the diagonal matrix with $4-k$ ones followed by $k$ $l$'s.
Suppose the eigenvalue of $T (l,k)$ on $\xi$ is $a(l,k)$, and define
the Hecke polynomial $H (\xi)$ of $\xi$ by
\[
H (\xi) = \sum_{k} (-1)^{k}l^{k (k-1)/2} a (l,k) T^{k} \in \C [T].
\]
Let $\rho \colon \Gal(\bar \Q / \Q)\rightarrow GL_{n} (\Q_{p})$ be a
continuous semisimple Galois representation unramified outside $pN$.
Then we say the eigenclass $\xi$ is \emph{attached} to $\rho$ if for
all $l$ not dividing $pN$ we have 
\[
H (\xi) = \det (1-\rho (\Frob_{l})T).
\]

The goal of this section is to formulate a conjecture about some
Eisenstein classes in the cohomology $H^{5} (\Gamma_{0} (N), \C)$.  We
recall the definition; for more background on Eisenstein cohomology
and its applications, we refer to \cite{harder-icm}.  

Let $X$ be the global symmetric space $SL (4,\R)/ SO (4)$, and let
$X^{\BorelSerre}$ be the partial compactification constructed by Borel
and Serre \cite{bs}.  The quotient $Y := \Gamma _{0}
(N)\backslash X$ is an orbifold, and the quotient $Y^{\BorelSerre } :=
\Gamma_{0} (N)\backslash X^{\BorelSerre}$ is a compact orbifold with
corners.  We have $H^{*} (\Gamma_{0} (N), \C) \simeq H^{*} (Y, \C)
\simeq H^{*} (Y^{\BorelSerre}, \C)$.

Let $\partial Y^{\BorelSerre} = Y^{\BorelSerre}\smallsetminus Y$.  The
Hecke operators act on the cohomology of the boundary $H^{*} (\partial
Y^{\BorelSerre}, \C)$, and the inclusion of the boundary $\iota \colon
\partial Y^{\BorelSerre} \rightarrow Y^{\BorelSerre}$ induces a map on
cohomology $\iota^{*}\colon H^{*} (Y^{\BorelSerre}, \C) \rightarrow
H^{*} (\partial Y^{\BorelSerre}, \C)$ compatible with the Hecke
action.  The kernel $H^{*}_{!}  (Y^{\BorelSerre}, \C)$ of $\iota^{*}$
is called the \emph{interior cohomology}; it contains the cohomology
with compact supports.  The goal of Eisenstein cohomology is to use
Eisenstein series and cohomology classes on the boundary to construct
a Hecke-equivariant section $s\colon H^{*} (\partial Y^{\BorelSerre},
\C) \rightarrow H^{*} (Y^{\BorelSerre}, \C)$ mapping onto a complement
$H^{*}_{\Eis} (Y^{\BorelSerre}, \C )$ of the interior cohomology in
the full cohomology.  We call classes in the image of $s$
\emph{Eisenstein classes}.  (In general, residues of Eisenstein series
can give interior, noncuspidal cohomology classes, with infinity type
a Speh representation, but as noted in \cite{AGM}, these do not
contribute to degree 5.)

To describe our conjectural Eisenstein classes, we give the Galois
representations we believe are attached to the classes along with the
corresponding Hecke polynomials.  In the following, $\varepsilon$
denotes the $p$-adic cyclotomic character, so that $\varepsilon (\Frob
_l) = l$ for any prime $l$ coprime to $p$.  We denote the trivial
representation by $i$.

\begin{itemize}
\item \emph{Weight two holomorphic modular forms:} Let
$\sigma_{2}$ be the Galois representation attached to a holomorphic
weight 2 newform $f$ of level $N$ with trivial Nebentypus.  Let
$\alpha$ be the eigenvalue of the classical Hecke operator $T_{l}$ on
$f$.  Let $\IIa (\sigma_{2})$ and $\IIb (\sigma_{2})$ be the 
Galois representations in the first two rows of Table \ref{grtab} (see
p.~\pageref{grtab}).

\item \emph{Weight four holomorphic modular forms:} Let
$\sigma_{4}$ be the Galois representation attached to a holomorphic
weight 4 newform $f$ of level $N$ with trivial Nebentypus.  Let
$\beta$ be the eigenvalue of the classical Hecke operator $T_{l}$ on
$f$.  Let $\IV(\sigma_{4})$ be the Galois representation in the third
row of Table \ref{grtab}.

\item \emph{Cuspidal cohomology classes from subgroups of $SL
(3,\Z)$:} Let $\tau$ be the Galois representation conjecturally
attached to a pair of nonselfdual cuspidal cohomology classes
$\eta,\eta ' \in H^{3} (\Gamma^{*}_{0} (N), \C)$, where $\Gamma^{*}_0
(N)\subset SL (3,\Z)$ is the congruence subgroup with bottom row
congruent to $(0,0,*)$ modulo $N$.  Let $\gamma$ be the eigenvalue of
the Hecke operator $T_{l,1}$ on $\eta $, and let $\gamma '$ be its
complex conjugate.  Let $\IIIa(\tau)$ and $\IIIb (\tau)$ be the Galois
representations in the last two rows of Table \ref{grtab}.

\end{itemize}

If $f$ is a weight $2$ or weight $4$ eigenform as above, we denote by
$d_{f}$ the degree of the extension of $\Q$ generated by the
eigenvalues of $f$.  Similarly for an eigenclass $\eta \in H^{3}
(\Gamma^{*}_{0} (N), \C)$ we write $d_{\eta}$ for the degree of the
field generated by the eigenvalues of $\eta $.  Also, the $L$-function
$\Lambda (f,s)$ of a holomorphic modular form of even weight $k$ and level
$N$ refers to
the function 
\[
\Lambda (f,s) = N^{s/2} (2\pi)^{-s}\Gamma (s)L (f,s),
\]
where $L (f,s)$ is the Dirichlet series for $f$.  
The function $\Lambda (f,s)$ satisfies the functional equation 
\[
\Lambda (f,s) = w(-1)^{k/2}\Lambda (f,k-s),
\]
where $w\in \{\pm 1 \}$ gives the action of the Fricke
involution.

\begin{conjecture}\label{eisenstein.conjecture}
Fix a positive prime $p$.  Then the cohomology group $H^{5} (\Gamma_{0}
(p), \C)$ contains the following Eisenstein subspaces:
\begin{enumerate}
\item For each weight two holomorphic newform $f$ of level $p$ with
associated Galois representation $\sigma_{2}$, two $d_{f}$-dimensional
subspaces, one attached to the Galois representation $\IIa
(\sigma_{2})$, and the other to the Galois representation $\IIb
(\sigma_{2})$.
\item For each weight four holomorphic newform $f$ of level $p$ with
associated Galois representation $\sigma_{4}$ such that the sign
$w$ of
the functional equation of the $L$-function $\Lambda (f,s)$ is
negative, a $d_{f}$-dimensional subspace attached to the Galois
representation $\IV (\sigma_{2})$.
\item For each pair of nonselfdual cuspidal cohomology class $\eta,
\eta ' \in H^{3} (\Gamma^{*}_{0} (p), \C)$, $\Gamma^{*}_0 (p)\subset
SL (3,\Z)$ with conjecturally associated Galois representation $\tau$,
two $d_{\eta }$-dimensional subspaces, one attached to the Galois
representation $\IIIa (\tau)$, and the other to the Galois
representation $\IIIb (\tau)$.
\end{enumerate}
\end{conjecture}

\begin{example}
Let $N=53$.  Then $N$ is prime, and is in fact the first level for
which $H^{3} (\Gamma^{*}_{0} (N), \C)$ contains nontrivial cuspidal
classes.  According to \cite{AGG}, there are two nonselfdual cuspidal
classes $\eta , \eta '$ whose Hecke eigenvalues are complex conjugates
of each other.  Moreover if $T (l,1)\eta = a (l,1)\eta$, then one
knows that $T (l,2)\eta = \bar a (l,1)\eta$.  Writing $\omega =
(1+\sqrt{-11})/2$, the Hecke eigenvalues of $\eta$ are given by the
following table:
\begin{center}
\begin{tabular*}{0.75\textwidth}{@{\extracolsep{\fill}}c||c|c|c|c|c|c}
$l$&2&3&5&7&11&13\\
\hline
$a (l,1)$&$-1-2\omega$&$-2+2\omega$&$1$&$-3$&$1$&$-2-12\omega$\\
\end{tabular*}
\end{center}
According to Conjecture \ref{eisenstein.conjecture}, this pair should
contribute a $4$-dimensional subspace to $H^{5} (\Gamma_{0} (53), \C)$.

Now we consider the subspaces corresponding to modular forms.  By
consulting tables of modular forms of weights 2 and 4 \cite{stein}, we
find (i) the dimension of the space of weight 2 newforms of level $53$
is $4$, with one form rational and one defined over the real cubic
field of discriminant $148$, and (ii) the dimension of the space of
weight 4 newforms of level $53$ is $13$, with one form rational, one
defined over a real quartic field, and one defined over a real octic
field.  We can also see that the rational and quartic weight 4 forms
have negative sign in their functional equations.  Hence the weight 2
forms should contribute a $2 (1+3)=8$ dimensional subspace, while the
weight 4 forms should contribute a $1+4=5$ dimensional subspace.

Thus the final result predicted by Conjecture
\ref{eisenstein.conjecture} is that $\dim H^{5} (\Gamma_{0} (53), \C)
\geq 4+8+5 = 17$.  Indeed, our computations show $\dim H^{5}
(\Gamma_{0} (53), \C) = 17$, and that the Hecke polynomials of the
eigenclasses match those predicted by Conjecture
\ref{eisenstein.conjecture} at $l=2$.  In particular, there is no
cuspidal cohomology at level $53$.

\end{example}

We remark that if we allow $N$ to be composite, we know that
Conjecture \ref{eisenstein.conjecture} does not give a complete
description of the Eisenstein subspace of $H^{5} (\Gamma_{0} (N),\C)$.
For instance, the cohomology contains Eisenstein classes corresponding
to minimal faces of the Borel--Serre boundary when $N$ is a square, as
mentioned on p.~\pageref{minimalref}.  Moreover, newforms for levels
properly dividing $N$ also appear in $H^{5}$.  As an example, it
appears from our data that if a prime $p$ exactly divides $N$, then a
weight $2$ newform $f$ at level $p$ contributes two
$3d_{f}$-dimensional subspaces to $H^{5} (\Gamma_{0} (N),\C)$,
corresponding to three copies of the representations $\IIa$ and
$\IIb$.  Similar phenomena occur for weight $4$ forms.  Also, for
sufficiently composite levels, Eisenstein cohomology involving weight
3 forms with odd character can occur, as happens e.g. for level 50.
Since we don't fully understand the mechanisms underlying these lifts,
we restrict Conjecture \ref{eisenstein.conjecture} to prime level.

\section{Cuspforms}\label{cusp}

Let $\pi$ be a cuspidal automorphic representation for $GL(4,\A)$,
where $\A$ denotes the adele group of $\Q$.  Assume that $\pi$
contributes to the cohomology $H^5(\Gamma_0(N),\C)$.  Then the
infinity type $\pi_\infty$ is uniquely determined, and is denoted
$\pi_1$ in the table of \cite[p.~65]{AG}.  As explained there, when
restricted to $SL(4,\R)$, $\pi_1$ breaks up into the direct sum
$\pi_1^+ \oplus \pi_1^-$, whose components are interchanged by the
inner automorphism $\iota$ induced by $\diag(-1,1,1,1)$.

It follows from this last fact that $H^5_{\cusp}(\Gamma_0(N),\C)$, the
cuspidal part of the cohomology, will have isotypic components for the
action of the Hecke algebra of even dimension $2k$, and $\iota$ will
act as an involution on each isotypic component interchanging two
complementary subspaces of dimension~$k$.  (Compare what happens for
$GL(2)$, where every cohomological cuspidal automorphic representation
contributes to the group cohomology twice, once as a holomorphic
modular form and once as an anti-holomorphic form.)

Let $f$ be a Hecke eigenclass in $H^5_{\cusp}(\Gamma_0(N),\C)$.  For a
fixed prime $\ell$ not dividing $N$, let $a,b,c$ be the Hecke
eigenvalues of $T_{\ell,1},T_{\ell,2},T_{\ell,3}$ respectively.  Then
the Hecke polynomial at $\ell$ is by definition
$$
P(X)=1-aX+b\ell X^2 -c\ell^3 X^3 + \ell^6 X^4.
$$

Suppose $\pi$ is the cuspidal automorphic representation associated to
$f$.  We refer to \cite[pp.~756--7]{APT} for the following facts.
Letting $c$ denote complex conjugation, there is defined another (or
possibly the same) cuspidal automorphic representation $\cpi$ with the
property that the Hecke eigenvalues of the corresponding cohomology
class ${}^c f$ are $\bar a, \bar b, \bar c$.

There is also the contragredient cuspidal automorphic representation
$\tpi$ and the Hecke eigenvalues of its corresponding cohomology class
$\widetilde f$ are $c,b,a$.  Because the coefficients of our
cohomology class are trivial, the weight $w$ in the notation of
\cite{APT} equals 3 and $\cpi=\tpi$.  Therefore $a=\bar c$ and $b =
\bar b \in \R$.  We say that $f$ or $\pi$ are \emph{selfdual} if $\tpi
\simeq \pi$.  This happens if and only if $a=c$ and hence if and only
if $a,b,c\in \R$.

One can recognize nonzero elements in $H^5_{\cusp}(\Gamma_0(N),\C)$ as
follows.  Compute Hecke eigenvalues on the whole cohomology
$H^5(\Gamma_0(N),\C)$.  Any system of eigenvalues that does not appear
to be attached to an Eisenstein cohomology class must be attached to
cuspidal cohomology.  For example, if even a single Hecke polynomial
is irreducible, then the corresponding Hecke eigenspace must be
cuspidal.  A further check is given by the fact that a cuspidal
eigenspace must be even dimensional.

Let $V$ be a minimal nontrivial Galois-stable Hecke eigenspace in
$H^5_{\cusp}(\Gamma_0(N),\overline{\Q})$.  The Hecke eigenvalues on
$V$ generate an order $R$ in the ring of integers in some number
field.  If $R$ is totally real then the corresponding automorphic
representations are self-dual and should be functorial liftings from a
smaller group that fixes a quadratic form, i.e. either from $GSp(4)$
or $GO(4)$.  Otherwise, $R$ must generate a complex CM field (see
\cite{APT}).

For more information about the situation where $R$ is totally real, in
which the desired results are close to being proved, we refer to
\cite{rsha, rproc}. In brief, assume $\pi$ is essentially selfdual,
i.e.~$\tilde \pi$ is isomorphic to $\pi \otimes \chi$ for some
character $\chi$. The $L$-groups of $GSp(4)$ and $GO(4)$ can be
identified, respectvely, with their complex points.  Then $\pi$ should
descend to a cusp form $\Pi$ on $GO(4)/\Q$ (resp. $GSp(4)/\Q$) iff the
symmetric square (resp.~the exterior square) L-function of $\pi$
admits, when twisted by the inverse of a character $\nu$ (``similitude
norm"), a pole at $s=1$. (This corresponds to the symmetric
(resp.~exterior) square of the associated 4-dimensional representation
$\sigma$ of the conjectural Langlands group $\L_\Q$ having a stable
line.)

We found two-dimensional spaces of cuspidal cohomology at levels $N =
61,73,79$, and in each case the Hecke polynomial at 2 was irreducible.
The rest of the cohomology at these levels is accounted for by the
Eisenstein subspaces of Conjecture \ref{eisenstein.conjecture}.  Since
$V$ in each of these cases is 2-dimensional, the Hecke eigenvalues in
each case must be rational integers.  Therefore these cuspforms are
selfdual and are expected to be lifts form $GSp(4)$ or $GO(4)$.

If we assume the Weil bounds for our Hecke eigenvalues, they tell us
that the eigenvalues $a,c$ at $\ell$ all have absolute value less than
or equal to $4\ell^{3/2}$ and $|b| \le 6\ell^2$.  Hence although we
only work modulo 31991 or 12379, we can assert the eigenvalues as
found in Table \ref{cusptab}.  For level $79$, we encountered overflow errors
working modulo 31991, and instead we worked modulo 12379.

If the Hecke eigenvalues were to generate for example an imaginary
quadratic extension of $\Q$, $V$ would have to be at least 4
dimensional and the corresponding cuspforms would not be lifts from
smaller groups.  It would be of great interest to find examples of
this.  The analogous objects do exist for $GL(3)$ as first found in
\cite{AGG}.

As stated above, the cuspforms we found are expected to be lifts from
$GSp(4)$ or $GO(4)$.  In fact, thanks to an argument shown us by
Ramakrishnan, we can show that for these levels, our cuspforms are
always lifts from $GSp(4)$:

\begin{theorem}\label{Dinakar}
Let $f$ be a Hecke eigenclass in $H^5_{\cusp}(\Gamma_0(N),\C)$ and
$\pi$ the cuspidal automorphic representation associated to
$f$. Assume that $\pi$ is a functorial lift from $GO(4)$.  Then $N$
cannot be squarefree.
\end{theorem}

\begin{proof} (D. Ramakrishnan)
The archimedean parameter of $\pi$ is a homomorphism
$$
\sigma_\infty: W_\R \to GL(4,\C),
$$
where $W_\R$ is the real Weil group containing $\C^*$ as a subgroup and
$\Gal(\C/\R)$ as the corresponding quotient, whose non-trivial element
$c$ acts on $\C^*$ by sending $z$ to its complex conjugate
$\overline{z}$. Since $\pi$ contributes to cuspidal cohomology with
constant coefficients of the congruence subgroup $\Gamma$ of
$SL(4,\Z)$, one necessarily has (in the unitary normalization):
$$
\sigma_\infty \simeq I(W_\R,\C^*;\alpha^3)\oplus
I(W_\R;\C^*;\alpha),
$$
where $I$ denotes induction, here from $\C^*$ to $W_\R$, and
$\alpha=z/|z|$.  Consequently the restriction of $\sigma_\infty$ to
$\C^*$ is the sum of the characters in the ``infinity type":

$$ p_\infty = \{\alpha^3,\alpha,\alpha^{-1},\alpha^{-3}\}.
$$

Suppose $\pi$ is of general orthogonal type. Then it is either of the
following two kinds:
$$ (I) \ \ \ \ \ \pi = \pi_1
\boxtimes \pi_2,
$$
where $\pi_1, \pi_2$ are cusp forms on $GL(2)/\Q$ and $\boxtimes$ is
the Rankin--Selberg (or automorphic) tensor product, which corresponds
to the tensor product (not the direct sum) of the corresponding
2-dimensional representations of $\L_\Q$; or
$$ (II) \
\ \ \ \ \pi = \As_{K/\Q}(\eta),
$$
where $\pi$ is the Asai representation defined by a cusp form $\eta$
of $GL(2)/K$, where $K$ is a quadratic extension of $\Q$.

In case (I), suppose one of the $\pi_j$, say $\pi_2$, is an Eisenstein
class of the form $\mu_1\boxplus \mu_2$ (``isobaric sum").  This means
$L(s,\pi_{2})=L(s,\mu_1)L(s,\mu_2)$, which implies $\pi =
(\pi_1\otimes\mu_1) \boxplus (\pi_1\otimes\mu_2)$, and thus $\pi$ is
certainly not cuspidal.

Continuing with case (I), assume now that both the $\pi_j$ are
cuspdial and let $\sigma_{j,\infty}$ denote the
$W_\R$-parameter of the cusp form $\pi_j$. Since the two irreducible
constituents of $\sigma_\infty$ are not twist equivalent,
$\sigma_{1,\infty}$ and $\sigma_{2,\infty}$ are both forced to be
irreducible. We may write, after possiby interchanging $\pi_1$ and
$\pi_2$,

$$ \sigma_{1,\infty}= I(W_\R,\C^*;
\alpha^a)
$$
and
$$
\sigma_{2,\infty}= I(W_\R,\C^*;
\alpha^b)
$$
with $a \geq b > 0$.

Since the Rankin--Selberg product $(\pi_1,\pi_2) \to \pi$ is functorial
at all places \cite{ramaannals}, in particular at
$\infty$, we must have
$$
\sigma_\infty =
\sigma_{1,\infty}\otimes\sigma_{2,\infty},
$$
implying that
$$
a+b=3, a-b=1,
$$
so that $a=2, b=1$. In other words, $\pi_1, \pi_2$ are classical
holomorphic newforms of weight $3, 2$ respectively. Let $N_j$ be the
level (conductor) of $\pi_j$ (for $j=1,2$) and $N$ the conductor of
$\pi$. We have
$$
(N_1,N_2)=1
\implies N=N_1^2N_2^2
$$
(see \cite{br}).  Since we are assuming $N$ is
squarefree, this cannot happen.

More generally, if at any prime $p$, $v_p(N_1)>0$ and $v_p(N_2)=0$,
then $v_p(N) = 2v_p(N_1)$. So for $N$ to be square-free, it is
necessary that $v_p(N_1)=v_p(N_2)$ at all $p$. In other words,
$v_p(N_1)$ and $v_p(N_2)$ are simultaneously 0 or simultaneously 1.

Now let $p$ divide $N$. It is then the case that up to
an unramified twist that can be ignored, $\pi_{1,p}\simeq\pi_{2,p}\simeq St_p$.
(See, for example, the table on p. 73 of \cite{gelbart} giving the conductors
of representations of $GL(2,\Q_p)$; note also that the conductor of
any ramified twist of $St_p$ is divisible by $p^2$.)  The reason is
that
$$
St_p \boxtimes St_p = sym^2(St_p) \boxplus
1,
$$
where $sym^2(St_p)$, being the Steinberg representation of
$GL(3,\Q_p)$, is also of conductor $p$. Indeed, $St_p$
corresponds, by the local correspondence at $p$, to
$$
\tau_p
=1\otimes id: W_{\Q_p}\times SL(2,\C) \to GL(2,\C), (w,g)\to
g,
$$
and the Steinberg representation of $GL(3,\Q_p)$ corresponds to
$$
1\otimes
sym^2: W_{\Q_p} \times SL(2,\C) \to GL(3,\C).
$$
The only other possible representations of $GL(2,\Q_p)$ of conductor
$p$ are the principal series representation $\xi_1 \boxplus \xi_2$,
with $\xi_1$ of conductor $p$ and $\xi_2$ unramified. But if
$\pi_{2,p}$ is of this form with $\pi_{1,p}=St_p$, the conductor of
$\pi_p$ will be divisible by $p^2$. Similarly, if $\pi_{1,p}$ and
$\pi_{2,p}$ are both principal series of conductor $p$, then the
conductor of $\pi_p$ will be divisible by $p^2$.

In our case, $\pi_1$ is generated by a holomorphic newform of weight
3, hence has a non-trivial character since its character must have the
same parity as the weight of $\pi_1$, and hence must be odd.  This
character must then be ramified at some prime $p_0$, say, because $\Q$
has class number 1. This $p_0$ must divide $N$, and we get a
contradiction from $\pi_{1,p_0}$ being $St_{p_0}$, which has trivial
central character.  Note that this argument depends on $N$ being squarefree.

Now we can move to case (II). Suppose $\pi=\As_{K/\Q}(\eta)$, for a
cusp form $\eta$ on $GL(2)/K$, $K$ a quadratic field. A basic property
of Asai representations gives \cite{ramaartin}
$$
\pi_K = \eta \boxtimes
(\eta^\theta),
$$
where $\pi_K$ denotes the base change of $\pi$ to $GL(4)/K$, $\theta$
is the nontrivial automorphism of $K$, and $\eta^\theta$ means
$\eta\circ\theta$.  Let $\tau_\infty$ denote the $2$-dimensional
representation of $W_{K_\infty}$ asociated to $\eta_\infty$, so that
$\tau_\infty^{\theta}$ is asociated to $\eta_\infty^\theta$.

First consider the case of a real quadratic K. Then we have
$$
\sigma_\infty = \tau_\infty \otimes
\tau_\infty^\theta.
$$
Arguing as in case (I) we see that $\tau_\infty$ should have parameter
$\{\alpha^2,\alpha^{-2}\}$, while $\tau^\theta_\infty$ has parameter
$\{\alpha,\alpha^{-1}\}$, which is clearly impossible.

So we may assume that $K$ is imaginary quadratic where $\theta$
induces complex conjugation on $K_\infty=\C$. Suppose the archimedean
parameter of $\eta$ is $\{\alpha^a, \alpha^c\}$.  Since $W_\C=\C^*$,
this need not be preserved by complex conjugaton. But nevertheless, it
forces the archimdean parameter of $\eta^\theta$ to be
$\{\alpha^{-a},\alpha^{-c}\}$, and the tensor product of these two
parameters is not regular, though tempered, and thus cannot contribute
to the cohomology of an arithmetic group.
\end{proof}

We remark that from \cite[p.~52]{AG} we see that the automorpic
representation on $GSp(4)$ that we lift to get our $f$ must correspond
to a holomorphic Siegel modular form of weight 3.

\newpage
\section*{Tables of Results}\label{table}

\begin{table}[ht]
\begin{center}
\begin{tabular}{|p{30pt}|p{30pt}||p{30pt}|p{30pt}||p{30pt}|p{30pt}|}
\hline
Level&rank&Level&rank&Level&rank\\
\hline\hline
2&0&23&5&44&18\\
3&0&24&2&45&27\\
4&0&25&7&46&19\\
5&0&26&7&47&11\\
6&0&27&12&48&26\\
7&0&28&7&49&20\\
8&0&29&6&50&34\\
9&3&30&8&51&19\\
10&0&31&6&52&21\\
11&2&32&12&53&17\\
12&0&33&10&54&49\\
13&1&34&12&55&15\\
14&2&35&7&56&20\\
15&2&36&24&57&19\\
16&3&37&8&59&14\\
17&3&38&14&61&20\\
18&9&39&10&67&17\\
19&3&40&9&71&17\\
20&2&41&9&73&20\\
21&3&42&17&79&25\\
22&7&43&10&83&21\\
\hline
\end{tabular}
\end{center}
\caption{Betti numbers for $H^{5} (\Gamma_{0} (N),
\C)$.\label{bettitab}}
\end{table}

\begin{table}[ht]
\begin{center}
\begin{tabular}{|p{40pt}||p{80pt}|p{150pt}|}
\hline
\IIa & $\sigma_{2}\oplus \varepsilon^{2}\oplus \varepsilon^{3}$ & $(1-l^{2}T) (1-l^{3}T) (1-\alpha T + lT^{2})$\\
\hline
\IIb & $i\oplus \varepsilon^{2}\sigma_{2}\oplus \varepsilon  $ &$(1-T) (1-l T)(1-l^{2}\alpha T + l^{5}T^{2})$\\
\hline
\hline
\IV & $\sigma_{4}\oplus \varepsilon \oplus \varepsilon^{2}$ & $(1-lT) (1-l^{2}T) (1-\beta T + l^{3}T^{2})$\\
\hline
\hline
\IIIa & $\tau \oplus \varepsilon^{3} $&$(1-l^{3}T) (1-\gamma T + l \gamma 'T^{2} - l^{3}T^{3})$\\
\IIIb & $i\oplus \varepsilon \tau$&$(1-T) (1-l\gamma T + l^{3} \gamma 'T^{2} - l^{6}T^{3})$\\
\hline
\end{tabular}
\end{center}
\caption{Galois representations and Hecke polynomials for Eisenstein
classes\label{grtab}.  See Section \ref{eis} for explanation of notation.}
\end{table}

\begin{table}[htb]
\begin{center}
\begin{tabular}{|p{40pt}||p{80pt}|p{150pt}|}
\hline
\multicolumn{3}{|l|}{Level \textbf{61}}\\
\hline
$T_{2}$&$(-7,12,-7)$&$1+7T+24T^{2}+56T^{3}+64T^{4}$\\
\hline
$T_{3}$&$(-3,1,-3)$&$1+3T+3T^{2}+81T^{3}+729T^{4}$\\
\hline
\hline
\multicolumn{3}{|l|}{Level \textbf{73}}\\
\hline
$T_{2}$&$(-6,11,-6)$&$1+6T+22T^{2}+48T^{3}+64T^{4}$\\
\hline
$T_{3}$&$(-2,1,-2)$&$1+2T +3T^{2}
+54T^{3} +729T^{4}$\\
\hline
\hline
\multicolumn{3}{|l|}{Level \textbf{79}}\\
\hline
$T_{2}$&$(-5,7,-5)$&$1+5T+14T^{2}+40T^{3}+64T^{4}$\\
\hline
$T_{3}$&$(-5,14,-5)$& $1+5T+42T^{2} + 135T^{3} + 729 T^{4}$\\
\hline
\end{tabular}
\end{center}
\caption{Eigenvalues and Hecke polynomials for $H^{5}_{\cusp}
(\Gamma_{0} (N), \C )$.\label{cusptab}}
\end{table}


\begin{thebibliography}{99}



\bibitem{AG} Avner Ash and David Ginzburg, \emph{$p$-adic
$L$-functions for $GL(2n)$}, Inv. Math. 116 (1994) 27--73.

\bibitem{AGG} Avner Ash, Daniel Grayson and Philip Green,
\emph{Computations of Cuspidal Cohomology of Congruence Subgroups of
$SL(3,\Z)$}, J. Number Theory 19 (1984) 412--436.

\bibitem{AGM} Avner Ash, Paul E. Gunnells and Mark McConnell,
\emph{Cohomology of Congruence Subgroups of $SL_4(\Z)$}, J. Number
Theory 94 (2002) 181--212.

\bibitem{APT} Avner Ash, Richard Pinch and Richard Taylor, \emph{An
$\hat A_4$ extension of $\Q$ attached to a non-selfdual automorphic
form on $GL(3)$}, Math. Ann. 291 (1991) 753--766.

\bibitem{br}
Laure Barthel and Dinakar Ramakrishnan,
\emph{A nonvanishing result for twists of $L$-functions of ${\rm GL}(n)$},
Duke Math. J. 74 (1994), no. 3, 681--700. 

\bibitem{bs} Armand Borel and J.-P. Serre, \emph{Corners and
arithmetic groups. Avec un appendice: Arrondissement des variétés à
coins, par A. Douady et L. Hérault}, Comment. Math. Helv. 48 (1973),
436--491.

\bibitem{vgt} Bert van Geemen and Jaap Top, \emph{A non-selfdual
automorphic representation of ${\rm GL}_3$ and a Galois
representation}, Invent. Math. 117 (1994), no. 3, 391--401.


\bibitem{gelbart} Stephen S. Gelbart, \emph{Automorphic forms on adele
groups}, Annals of Mathematics Studies, No. 83.  Princeton University
Press, Princeton, N.J.

\bibitem{grit} Valeri A. Gritsenko, \emph{Irrationality of the moduli spaces of
polarized abelian surfaces}, Proc. Egloffstein conf. on Abelian
Varieties, de Gruyter, Berlin (1995), 63--81.


\bibitem{harder-icm}
G\"unter Harder, \emph{Eisenstein cohomology of arithmetic groups and its applications
  to number theory}, Proceedings of the International Congress of
  Mathematicians, Kyoto, 1990 (Tokyo), Math. Soc. Japan, 1991,
  779--790.

\bibitem{ibuk} Tomoyoshi Ibukiyama, \emph{Siegel modular forms of
weight three and conjectural correspondence of Shimura type and
Langlands type}, in The conference on L-functions, Ed. L. Weng and
M. Kaneko, World Scientific (2007), 55--69.

\bibitem{PY} Cris Poor and David S. Yuen, \emph{Dimensions of cusp
forms for $\Gamma_0(p)$ in degree two and small weights}, preprint.

\bibitem{rproc} Dipendra Prasad and Dinakar Ramakrishnan, \emph{On the
global root numbers of $GL(n) \times GL(m)$}, in Automorphic Forms,
Automorphic Representations and Arithmetic, Proc. Symp. Pure Math. 66
part 2 (1999), 311--330.

\bibitem{ramaartin} Dinakar Ramakrishnan, \emph{Modularity of solvable
{A}rtin representations of {${\rm GO}(4)$}-type},
Int. Math. Res. Not. (2002), no. 1, 1--54.

\bibitem{ramaannals} Dinakar Ramakrishnan, \emph{Modularity of the
Rankin--Selberg $L$-series, and multiplicity one for ${\rm SL}(2)$},
Ann. of Math. (2) 152 (2000), no. 1, 45--111.

\bibitem{rsha} Dinakar Ramakrishnan and Freydoon Shahidi, \emph{Siegel
modular forms of genus 2 attached to elliptic curves}, to appear in
Mathematical Research Letters.

\bibitem{stein} William A. Stein, Tables of modular forms, available
from \\
\texttt{modular.math.washington.edu/Tables/tables.html}.

\end{thebibliography}
\end{document}